\documentclass[12pt,a4paper]{article}
\usepackage{amsmath,amsfonts,amssymb,amsthm}

\input epsf

 \small

\newcommand{\Af}{\mathrm{Af\mbox{}f}}
\newcommand{\af}{\mathrm{af\mbox{}f}}

\newcommand{\Span}{\mathrm{span}} 
\newcommand{\R}{\mathbb{R}} 
\newcommand{\Rf}{\R^4}

\newcommand{\ba}{\begin{array}}
\newcommand{\ea}{\end{array}}
\newcommand{\G}{Gl(4,\R)}
\newcommand{\Aff}{\Af(\Rf)}
\newcommand{\aff}{\af(\Rf)}
\newcommand{\spafo}[4]{\left(\! \ba{c} #1 \\ #2 \\ #3 \\ #4 \ea\! \right)}
\newcommand{\spafi}[5]{\left(\! \ba{c} #1 \\ #2 \\ #3 \\ #4 \\ #5 \ea\!
  \right)}
\newcommand{\vpi}{\varpi}
\newcommand{\F}{{\cal F}}
\newcommand{\iz}{\varphi}
\newcommand{\ic}[3]{\iz_{#1 #2}^{#3}}
\newcommand{\om}{\omega}
\newcommand{\np}{\nabla^\bot}
\newcommand{\pu}{\frac{\partial}{\partial u}}
\newcommand{\pv}{\frac{\partial}{\partial v}}
\newcommand{\tS}{\tilde{S}}
\newcommand{\tX}{\tilde{X}}
\newcommand{\tU}{\tilde{U}}
\newcommand{\tu}{\tilde{u}}
\newcommand{\tv}{\tilde{v}}
\newcommand{\tx}{\tilde{x}}
\newcommand{\tih}{\tilde{h}}
\newcommand{\tal}{\tilde{\alpha}}
\newcommand{\tbe}{\tilde{\beta}}
\newcommand{\tk}{\tilde{k}}
\newcommand{\al}{\alpha}
\newcommand{\be}{\beta}
\newcommand{\ga}{\gamma}
\newcommand{\sig}{\sigma}
\newcommand{\eps}{\varepsilon}
\newcommand{\tran}{\,^{T}\!}
\newcommand{\half}{{\tfrac{1}{2}}}
\newcommand{\scs}{semiconformal structure }
\newcommand{\Sym}{Sym(2)}

\newcommand{\spl}{Sl(2,\R)}
\newcommand{\iso}{O_{1}(3)}
\newcommand{\isop}{O_{1}^{++}(3)}
\newcommand{\Min}{\R^3_1}
\newcommand{\mat}[4]{\left(\! \ba{cc} #1 & #2 \\ #3 & #4 \ea\!
  \right)}
\newcommand{\pgl}{PGl^{+}(2,\R)}
\newcommand{\Id}{I\!d}

\newtheorem{lem}{Lemma}
\newtheorem*{thm}{Theorem}

\theoremstyle{remark}
\newtheorem*{rem}{Remark}

\newcounter{rom}
\renewcommand{\therom}{(\roman{rom})}
{\end{list}}
\begin{document}
\title{Parallel surfaces in affine 4-space}
\author{C. Scharlach\thanks{Partially supported by the DFG-project SI
163-7} \and L. Vrancken\thanks{Partially supported by a research
fellowship of the Alexander von Humboldt Stiftung (Germany)}} \maketitle

\begin{abstract}
We study affine immersions as introduced by Nomizu and Pinkall. We
classify those affine immersions of a surface in $\R^4$ which are
degenerate and have vanishing cubic form (i.e.  parallel second
fundamental form). This completes the classification of parallel
surfaces of which the first results were obtained in the beginning of
this century by Blaschke and his collaborators.
\end{abstract}

\medskip\noindent 
{\bfseries Subject class: } 53A15 
 

\section{Introduction}\label{sec:intro}
We consider the standard affine space $\R^{m}$ equipped with its
standard connection $D$.  Let $M^n$ be a manifold equipped with a
torsion free affine connection $\nabla$ and let $x:
(M^n,\nabla)\to (\R^{m},D)$, $m \geq n$ be an
immersion. Following \cite{nopi87}, we call $x$ an affine immersion if
there exists a transversal $(m-n)$-dimensional bundle $\sig$ such that
\begin{equation}
D_X x_*(Y) -x_*(\nabla_X Y) \in \sigma,\label{eq1}
\end{equation}
for all vector fields $X$ and $Y$ which are tangent to $M^n$. It is
immediately clear that if we equip $\mathbb R^{m}$ with a
semi-Riemannian metric and take for $\sigma$ the normal bundle, then
isometric immersions provide examples of affine
immersions. Also the equiaffine immersions, in the sense of Blaschke
for hypersurfaces, and in the sense of \cite{we37} or
\cite{NV} for higher codimensions provide examples of affine
immersions.

For an affine immersion it is possible to introduce a bilinear form
$h$, called the second fundamental form, which takes values in the
tranversal bundle $\sigma$ by
\begin{equation}
h(X,Y)=D_X x_*(Y) -x_*(\nabla_X Y) \in \sigma.\label{eq2}
\end{equation}
Since $\nabla$ is a torsion free affine connection, $h$ is symmetric in
$X$ and $Y$. Let $\xi$ be a vector field which takes values in
$\sigma$. Similarly, as for isometric immersions, we can now introduce
a normal connection $\nabla^\perp$ and Weingarten operators $A_\xi$ by
decomposing $D_X \xi$ into a tangential part and a part in the
direction of $\sigma$, i.e. we have the Weingarten formula which
states that
\begin{equation}
D_X \xi =-x_*(A_\xi X) + \nabla^\perp_X \xi.\label{equ4}
\end{equation}  

Using the Weingarten formula, it is now possible to define the
covariant derivative $\nabla h$ of the second fundamental form $h$ by
\begin{equation}
(\nabla_X h)(Y,Z)= \nabla_X^\perp h(Y,Z) - h(\nabla_X Y,Z) -
h(Y,\nabla_XZ).\label{equ3}
\end{equation}
Affine immersions for which $\nabla h$ vanishes identically are called
parallel immersions. In Riemannian geometry, these immersions and
their generalisations have been studied by many people, an overview
can be found in \cite{Lu}. A general classification of the Euclidean
parallel submanifolds was obtained in \cite{Fe74}. As far as
we know it is still an open problem to classify the 
semi-Euclidean parallel submanifolds.

In this paper we will focus on surfaces, i.~e. the dimension of $M^n$
equals two. All results will be local and valid on a suitable open
dense subset of $M^2$. We say that an affine immersion is linearly
full provided that for every point $p$ of $M^2$ and for every
neighborhood $U$ of $p$, $x(U)$ is not contained in a lower
dimensional affine subspace of $\mathbb R^n$. Using Lemma 2 of
\cite{Vr99} which says that $\sig=\text{im} h$ if a parallel affine
immersion is linearly full, it follows easily that a parallel surface
immersion which is linearly full has to be in $\mathbb R^2$, $\mathbb
R^3$, $\mathbb R^4$ or $\mathbb R^5$. The first case ($m=2$) clearly
implies that $M^2$ is an affine plane. In the other cases, a
nondegeneracy condition can be introduced as follows. Let $u=\{X_1,
X_2\}$ be a local basis in a neighborhood of a point $p$. Then we
define for $m=3$:
\begin{equation}
h_u(X,Y)=\det (X_1, X_2, D_X Y),
\end{equation}
and for $m=4$:
\begin{equation}
\label{confstr}
h_u(X,Y)=\half\left(\det (X_1, X_2, D_{X_1} X, D_{X_2} Y) + \det (X_1, X_2,
D_{X_1} Y, D_{X_2} X)\right).
\end{equation}
It is well known that in both cases the rank of $h_u$ is independent
of the choice of the local basis $u$. We call $M^2$ nondegenerate if
the rank equals 2, 1-degenerate if it equals 1 and 0-degenerate if it
equals 0. A surface in $\R^5$ is called nondegenerate if $\det (X_1,
X_2, D_{X_1} X_1, D_{X_1} X_2, D_{X_2} X_2) \neq 0$, which again is
independent of the choice of basis $u$. Using Lemma 2 of \cite{Vr99}
again, we see that a linearly full affine immersion of $M^2$ in
$\R^5$ is always nondegenerate.

Nondegenerate parallel immersions of a surface in $\mathbb R^3$,
$\mathbb R^4$ and $\mathbb R^5$ are considered in respectively
\cite{nopi87}, \cite{NV} and \cite{vramag}. Therefore, restricting to
an open and dense subset if necessary, only the degenerate cases still
need considering. If $M^2$ is contained in $\mathbb R^3$, a solution
was found in \cite{divr}. This leaves only the case that $M^2$ is
linearly full in $\mathbb R^4$ and degenerate. If $M^2$ is
0-degenerate and parallel, the immersion can not be linearly
full. Thus we are left with the 1-degenerate parallel surfaces which
are linearly full in $\R^4$. We prove the following:

\begin{thm}
Every 1-degenerate parallel affine surface immersion $(x,\sig)$ in
$\Rf$ is a ruled surface and can be locally parametrized either by
\begin{itemize} 
\item[{\rm I.1.}] $x(u,v)= \ga'(u) + v \ga(u)$, and 
  \newline $\sig=\Span(\xi_1,\xi_2)$ is
  given by \eqref{Fxi1} and \eqref{Fxi2}, or
\item[{\rm I.2.}] $x(u,v)= (\eps \gamma(u) + \gamma''(u)) + v \gamma'(u)$,
  $\eps=\pm 1$, and 
  \newline $\sig=\Span(\xi_1,\xi_2)$ is given by \eqref{Sxi1} and
  \eqref{Sxi2},
\item[{\rm II.}] $x(u,v)=\al(u)+v \be(u)$, $\be''=-\be$, and
  \newline $\sig=\Span(\xi_1,\xi_2)$ is given by \eqref{Txi1} and
  \eqref{Txi2}. 
\end{itemize}
\end{thm}

The paper is organized in two parts. In Section~\ref{sec:firstorder}
we apply the method of moving frames due to E. Cartan to an affine
immersion of $M^2$ in $\Rf$. We introduce the affine semiconformal
structure (cp.~\ref{confstr}), which was known already to
\cite{BM}. We end up with a classification of affine surfaces in $\Rf$
with respect to the (non)degeneracy-type of the affine semiconformal
structure and normal forms of the second fundamental form $h$ for each
type. This part is closely related to \cite{diss} and \cite{cafoinv}. 

In Section~\ref{sec:1-degPar} we restrict to 1-degenerate parallel
affine immersions of $M^2$ in $\Rf$. It turns out that they are ruled
surfaces (Lemma~\ref{ruled}) and therefore can be parametrized as
$x(u,v)= \al(u)+ v \be(u)$. We find special frames which simplify the
structure equations significantly. A reparametrization finally leads to
our main result.

We will use the Euler summation convention. 

\section{Classification of surfaces in $\Rf$ with respect to their affine
semiconformal structure}
\label{sec:firstorder}

\subsection{The affine frame bundle on $\Rf$}
\label{sec:affineframes}

We define a {\em frame} on $\Rf$ to be an ordered set 
\[ S_b=\{v_1,v_2,v_3,v_4;b\}, \quad\text{with}\quad
\spafo{v_1}{v_2}{v_3}{v_4} \in \G, b\in \Rf.\]

Let $F$ denote the set of all frames on $\Rf$ and $\pi :F \to \Rf$
the projection map, defined by: $\pi(S_b)=b$.
Let $\Aff$ be the Lie group of affine transformations on $\Rf$, 
\[ \Aff= \left\{ \left( \begin{tabular}{c|c} $A$&$0$\\
\cline{1-2} $b$&$1$\\ \end{tabular} \right)
: A\in \G, b\in \Rf\right\}. \]
Obviously we can identify $F$ with $\Aff$.
The local structure of $\Aff$ is encoded in the Lie algebra-valued
Maurer-Cartan form $\vpi_{|S}=dSS^{-1}\in \aff$, we use the notation:
\begin{equation}
  \label{strucEq}
  d\spafi{v_1}{v_2}{v_3}{v_4}{b}=
  \left( \begin{tabular}{c|c} $M$&$0$\\
      \cline{1-2} $n$&$0$\\ \end{tabular} \right) 
  \spafi{v_1}{v_2}{v_3}{v_4}{b}, \qquad M\in M(4\times 4,\R), n\in \Rf.
\end{equation}
If we let $\Aff$ act both on $F$ and $\{(b,1)\in \R^5 : b\in \Rf\}
\cong \Rf$ by {\em right multiplication} $R_S$, $R_S(C)=CS$, then $\pi
\circ R_S =R_S \circ \pi$.  To obtain the fibers of the bundle
$\F:=\pi \colon F \to \Rf$, note that $\pi(S_b)= (0,0,0,0,1)S_b$, and
the isotropy group of $(0,0,0,0,1)$ is
\[ H= \left\{ \left(  \begin{array}{c|c} A&0\\
\hline 0&1\\ \end{array} \right) \right\} \subset \Aff. \] 
We can identify the homogeneous space $\Aff/H \cong \Rf$ and obtain
that $\F$ is a principal (right) $H$-bundle, the {\em affine frame
  bundle} on $\Rf$.

\subsection{Adaption of the affine frame bundle to an affine surface immersion}
\label{sec:firstorderframe}

Let $U$ be a connected open subset of a two-dimensional oriented
manifold $M^2$ equipped with a torsion free affine connection $\nabla$
and let $x\colon (U, \nabla) \to (\Rf,D)$ be a smooth affine immersion
with transversal bundle $\sig$ (cp. Sec.~\ref{sec:intro}).  We want to
adapt the affine frame bundle to the surface by restricting the base
manifold to $x(U)$. We define the principal (right) $H$-bundle
$\F^0=\pi_U \colon F^0 \to U$ as the bundle over $U$ induced by $x$
and the affine frame bundle (cp. \cite[vol. V, p.  391f.]{Spivak} for
the notion of an induced bundle), i.~e. $\F^0= x^{*}\F$. For the
further adaption we take into account the given transversal bundle
$\sig$ and we use the first order information given by the tangent
bundle of $x$. We call a frame $S_u\in F^0$ a {\em first order frame}
if $\Span(v_1,v_2) = x_{*}(T_u M)$ and $\Span(v_3,v_4) = \sig$. We
denote the set of all first order frames on $U$ by $F^1\subset F^0$
and use the notation $S_u=\{v_1, v_2,\xi_1, \xi_2, x(u)\}\in F^1$.
The subgroup
\begin{equation}
  \label{H1}
  H^1=\left\{ \left( \ba{c|c}
    \ba{c|c} P & 0\\ \hline 0 & Q\\ \ea & 0 \\ \hline $0$ &$1$ \\ 
    \ea \right) : \det P \neq 0 \right\} \subset H, 
\end{equation}
acts transitively and effectively on $F^1$.  Thus we get a subbundle
$\F^1$ of $\F^0$, $\F^1=\pi_U\colon F^1 \to U$, where we use the same
notation for the restriction of $\pi_U$ to $F^1$. Obviously $\F^1$ is
a principal (right) $H^1$-bundle, the reduced bundle obtained by
reduction of the structure group $H$ of $\F^0$ to $H^1$ (cp.
\cite[vol.~I, pg.~53]{K-N}). Since the first two legs $v_1$ and $v_2$
of a frame $S_u \in F^1$ span the tangent space $x_{*}(T_u M)$, we get
two zero's in the last row of the Maurer-Cartan form $\vpi_{|{S_u}}$
on $F^1$ ($\vpi_5^3 =0$, $\vpi_5^4 =0$) and the forms $\vpi_5^1$ and
$\vpi_5^2$ drop down to $U$ (cp.~\ref{strucEq}). We use the notation:

\begin{equation}
  \label{strucEqFone}
  d\spafi{v_1}{v_2}{\xi_1}{\xi_2}{x}=
\left( \ba{c|c}
% Achtung, ziemliche Bastelei!!
    \ba{c|c} \ \iz\quad &\quad\psi \\ \hline \ \sigma\quad &\quad\tau \\
    \ea & 0 \\  
      \hline \ba{cccc} \!\!\om^1\!\!& \om^2\!&\  0& 0 \\ \ea & 0 \\ 
    \ea \right)
  \spafi{v_1}{v_2}{\xi_1}{\xi_2}{x}.
\end{equation}

For a fixed first order frame field (a smooth cross section of $\F^1$)
$S= \{v_1, v_2, \xi_1,\xi_2, x\}$, $v_i=dx(X_i)$, $X_i\in \Gamma(TM)$,
the entries of the Maurer-Cartan form define (resp. correspond to) the
following quantities ($1\leq i,j,k \leq 2$, $3 \leq \alpha \leq 4$)
(cp. \eqref{eq1}, \eqref{eq2} and \eqref{equ4}):
\begin{align}
  \nabla_{X_i} X_j &= \iz_j^k(X_i) dx(X_k) & \text{induced connection}
  \\
  h^\alpha(X_i,X_j) &= \psi_i^\alpha (X_j) & \text{second
    fundamental forms} \\
  -dx(A_{\xi_j}(X_i)) &= \sigma_j^k (X_i) dx(X_k) & \text{Weingarten
    operators} \\
  \np_{X_i} \xi_j &= \tau_j^\alpha(X_i) \xi_{(\alpha -2)} &
  \text{normal connection}
\end{align}
It is straightforward to show that $\np$ is a torsion-free affine
connection, $h^3$ and $h^4$ are symmetric bilinear forms and
$A_{\xi_1}$ and $A_{\xi_2}$ are 1-1 tensor fields.

\subsection{The affine semiconformal structure}
\label{sec:semiconformal}

To find the invariants (quantities independent of the choice of frame)
we can compute the (infinitesimal) group action (change of frames)
either on the Lie group level or on the Lie algebra level (cp.
\cite{G_equiv}) for the general theory, \cite{diss} for the
centroaffine case). Since the group of centroaffine transformations is
a subgroup of the affine group, the affine invariants are part of the
centroaffine ones. A description in detail can be found in
\cite{cafoinv}. We only will need the action on $\psi$ resp.  the
second fundamental forms $h^3$ and $h^4$. Let $S_u, \tS_u \in F^1$,
then there exists $B\in H^1$ such that $S_u= B\tS_u$ (cp.
\eqref{H1}). For the Maurer Cartan form we get:
\[ \vpi_{|S}B=dSS^{-1}B=d(B\tS)(B\tS)^{-1}B
=dB+Bd\tS\tS^{-1}=dB+B\vpi_{|\tS}.\] 
An evaluation of this equation and
$\psi_i^{\alpha}=h_{ij}^{\alpha}\om^j$ gives for $B= \left( \ba{c|c}
    \ba{c|c} P & 0\\ \hline 0 & Q\\ \ea & 0 \\ \hline $0$ &$1$ \\ 
    \ea \right)\in H^1$:
\begin{align}
\psi &= P \tilde{\psi} Q^{-1}, \label{cofh}\\
(h^3, h^4)&= (P[(Q^{-1})_1^1 \tih^3 + (Q^{-1})_2^1 \tih^4] \tran P,
P[(Q^{-1})_1^2 \tih^3 + (Q^{-1})_2^2 \tih^4]\tran P).\label{cof2ff}
\end{align}
For a frame $S_u=\{v_1,v_2,\xi_1,\xi_2,u\}\in F^1$ we define a
symmetric\footnote{We denote by $\odot$ the symmetric product of
  1-forms: $ \om \odot \eta(X,Y)=\half(\om(X)\eta(Y)+\om(Y)\eta(X))$.}
bilinear form $\phi$ on $F^1$ by:
\begin{equation} \label{defphi} 
\phi=\det \psi = \psi_1^3
\odot \psi_2^4 - \psi_2^3 \odot \psi_1^4 ,
\end{equation}
i.~e. $\phi(X,Y)=\frac12\left(\frac{[v_1,v_2,D_X dx(X_1),D_Y dx(X_2)]+
  [v_1,v_2,D_X dx(Y_1),D_X dx(X_2)])}{[v_1,v_2,\xi_1,\xi_2]}\right)$ for
some determinant form $[\quad]$ on $\Rf$ (cp.~\eqref{confstr}). 
We can use (\ref{cofh}) to determine how $\phi$ varies along the
fibers of $\F^1$: \begin{equation*} \label{cofscs} \phi = \det \psi =
  (\det P)(\det\tilde{\psi})(\det Q^{-1}) = \frac{\det P}{\det Q}
  \tilde{\phi}\: .
\end{equation*}
Now a {\em semiconformal structure} compatible with a quadratic
  form $q$ is defined as the set $\{rq : r\in
\R\!\setminus\!\!\{0\}\}$ and it makes sense to talk about a \scs
being nondegenerate, definite, etc.  (cp.~\cite[p. 4]{Weiner}).  The
quadratic form associated to $\phi$ induces a \scs on the tangent
space at each point of $U$. This structure on $U$ is called the {\em
  affine semiconformal structure} induced by $x$ and was known already
to \cite[p.~375]{BM}. Depending on the affine \scs we will call a
surface $x(U)$ a {\em nondegenerate, definite, indefinite} or {\em
 1-degenerate surface} if the induced affine \scs is nondegenerate,
definite, indefinite or 1-degenerate. A {\em 0-degenerate surface} is
a surface $x(U)$ for which the affine semiconformal structure contains
only the zero form.

\subsection{Normalization of $\psi$ and classification}
\label{sec:normalization2ff}

We saw that a change of frames induces an action of $H^1$
on $\Sym\!\times\!\Sym$ (cp. \eqref{cof2ff}), where $\Sym$ denotes the
algebra of all symmetric $2\!\times\!2$-matrices:
\begin{equation*}   \rho(B) (h^3,h^4) 
=(P[(Q^{-1})_1^1 h^3 + (Q^{-1})_2^1 h^4] \tran P,
P[(Q^{-1})_1^2 h^3 + (Q^{-1})_2^2 h^4]\tran P).
\end{equation*}

Note that the action can be written as the
composition of two actions of the form: 
\begin{align}
\label{action1} \rho_1(P) (h^3,h^4)&=(Ph^3\tran P,Ph^4\tran P), \\
\label{action2} \rho_2(Q)(h^3,h^4)&= ((Q^{-1})_1^1 h^3+(Q^{-1})_2^1 h^4,
(Q^{-1})_1^2 h^3+(Q^{-1})_2^2 h^4), \end{align} 
namely:
\begin{equation} \label{action}
\rho(\left( \ba{c|c}
    \ba{c|c} P & 0\\ \hline 0 & Q\\ \ea & 0 \\ \hline $0$ &$1$ \\ 
    \ea \right))=\rho_1(P)\circ\rho_2(Q).
\end{equation}

We want to choose normal forms (representatives of the orbits) in
$\Sym\!\times\!\Sym$ under the action of $H^1$ given by
(\ref{action}). Since the centroaffine situation is very close to the
affine one, we will omit some details. A more comprehensive
description can be found in \cite{cafoinv} (Section 4.1, 4.2). As we
just saw the action splits in two parts where $\Span (h^3,h^4)$ is an
invariant of the second part (\ref{action2}).  Therefore we want to
investigate the orbits of two-pencils\footnote{The span of two
  symmetric bilinear forms is called a two-pencil. } under the first
part (\ref{action1}) of the action. A first step is to understand the
action on a single element $h \in \Sym$. 

If we restrict $\rho_1$ to $\spl$, we can define an invariant quadratic
form $q$ in $\Sym$ by 
\begin{equation}\label{defq} q(h)=-\det h .\end{equation}
 Then $\Sym$ with the associated scalar product is isometric to the
Minkowski 3-space $\Min$ (see figure~\ref{lightcone}, for notations
cp. \cite{O'Neill}). This is easy to see if we choose
\begin{equation} \label{basis} E_0=\mat{1}{0}{0}{1},\: E_1= 
\mat{1}{0}{0}{-1},\: E_2=\mat{0}{1}{1}{0} \end{equation}
as a basis of $\Sym$. Then we get for every $h=aE_0+bE_1+cE_2=$
\begin{math}
  \bigl( \begin{smallmatrix}
        a+b&c\\c&a-b
        \end{smallmatrix} \bigr)
\end{math}
$\in \Sym$: 
$q(h)=-\det h=-(a^2-b^2-c^2)=-a^2+b^2+c^2$.

%%%%%%%%%%%%%%%%%%%%%%%%%%%%%%%%%%%%%%%%%%%%%%%%%%%%%%%%%%%%%%%%%
%%  Figures: \Figure{figure}{name}\Label\X, special to include ps-files:
%%    \psfile(llx,lly)(hsize,vsize)[scale,rotation]{file};
%%%%%%%%%%%%%%%%%%%%%%%%%%%%%%%%%%%%%%%%%%%%%%%%%%%%%%%%%%%%%%%%%
  \newdimen\psunit \psunit=1bp
\newcount\pshoffset \pshoffset=1 \newcount\psvoffset \psvoffset=1
\def\psfileA(#1,#2)(#3,#4)[#5,#6]#7{\immediate\write99{}%
  \immediate\write99{psfile data, #7: #5/100 (scale factor)}%
  \immediate\write99{\the\hsize (hsize,text), \the\vsize (vsize,text)}%
  \psunit=#3bp\message{#3bp=\the\psunit (hsize,figure),}%
  \psunit=#4bp\message{#4bp=\the\psunit (vsize,figure)}%
   \immediate\write99{}%
   \pshoffset=#1\multiply\pshoffset by #5\divide\pshoffset by 100%
   \psvoffset=#2\multiply\psvoffset by #5\divide\psvoffset by 100%
        \vbox to #4bp{\vss\hbox to #3bp{\kern0pt\relax%
        \includegraphics{#7}\hfil}}\ignorespaces}

\begin{figure}
\centerline{\framebox{\psfileA(55,580)(360,180)[200,0]{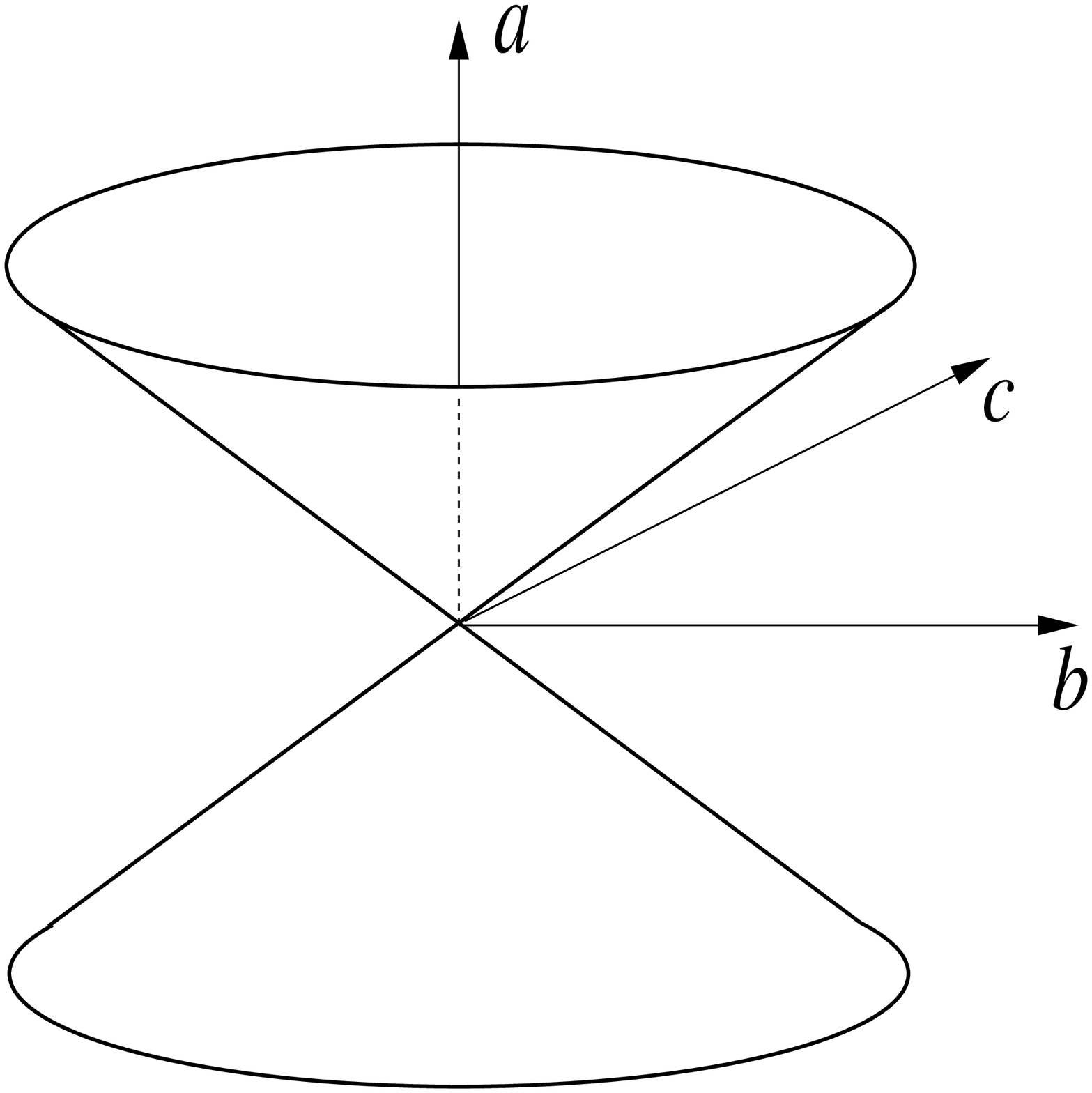}}}

\[ \ba{ccccc} 
h\: \hbox{spacelike} & \iff & \det h <0 &\iff&h\:\hbox{indefinite,} \\
h\: \hbox{timelike} & \iff & \det h >0 &\iff&h\:\hbox{definite,} \\
h\: \hbox{lightlike} & \iff & \det h =0 &\iff&h\:\hbox{degenerate.} 
\ea \]
\vspace{-3ex}
\caption{$(\protect\Sym,q)\protect\cong\protect\Min$}
\label{lightcone} 
\end{figure} 

Under this identification $\rho_1$ defines a representation of $\spl$ on
$\Min$. The invariance of $q$ means that $\rho_1 (P)$ is a linear isometry
of $\Min$, i.~e.  $\rho_1\colon \spl \to \iso$. This map is neither 1:1
($\rho_1 (P)=\rho_1(-P)$) nor onto ($\spl$ is connected, $\iso$ has four
components). However, $\rho_1\colon \pgl \to \R^+ \isop$ is an
isomorphism: 

\begin{thm}[\cite{cafoinv}, Thm.~3]
  Let $\isop$ be the group of all time- and space-orientation
  preserving isometries of $\Min$ and $\pgl =Gl^+(2,\R) /\{\pm\Id\}$,
  $Gl^+(2,\R)=\{P\in Gl(2,\R) : {\det P>0}\}$.  Identify $\Min$ and
  $\Sym$ by $\tran (a,b,c) \mapsto \mat{a+b}{c}{c}{a-b}$. Then
  $\rho_1\colon\pgl \to \R^+ \isop=\{rQ:\mbox{$r\in \R^+$}, Q\in
  \isop\}$, defined by $\rho_1 (P)A =PA\tran P$, is a (Lie group)
  isomorphism.
\end{thm}

Hence we know that $\pgl$ acts on an element $\bigl(
\begin{smallmatrix} a+b&c\\c&a-b \end{smallmatrix} \bigr) \in \Sym$ in
the same way as $\R^+ \isop$ acts on an element $\tran (a,b,c) \in
\Min$. The later action is well understood. $\isop$ acts transitively
on (ordered) orthonormal bases which have the same time- and 
space-orientation. Thus it acts also transitively on two-dimensional 
space-, time- or lightlike subspaces.

Two-pencils are either space-, time- or lightlike subspaces, either
two-dimensional or one-dim\-en\-sio\-nal or just the origin, where
dimension and type are invariant under $\rho_1$. Normal forms for the
lower dimensional cases are obvious \cite[p.~251]{Greub}.  In the
two-dimensional case we can choose the following normal
forms (cp.  (\ref{basis})):
\begin{align*}
\text{I.} && \Span (h^3,h^4)& \text{ spacelike:} && \Span (E_1,E_2), \\
\text{II.} && \Span (h^3,h^4)& \text{ lightlike:}&& \Span
(E_2,  \half(E_0+E_1)),  \\
\text{III.} && \Span (h^3,h^4)& \text{ timelike:} && \Span (E_0,E_1).
\end{align*}

Finally we can still use the second part $\rho_2$~\eqref{action2} of
the action $\rho$~\eqref{action} to map $h^3$ and $h^4$ to the
suitable basis vectors.

Summarized we obtain the following classes of surfaces in $\Rf$: 

\begin{minipage}{\textwidth}\label{surfacetypes}
\renewcommand{\arraystretch}{1.5}
\begin{center}
\begin{tabular}{rr|c|c|c|} 
& & $\Span(h^3,h^4)$ & $\phi$ & normal form \\ \hline
I. & & spacelike plane & definite & $(E_1,E_2)$ \\ \hline
II. & & lightlike plane & 1-degenerate & $(E_2, \frac12 (E_0+E_1))$ \\ \hline 
III. & & timelike plane & indefinite & $(E_0, E_1)$ \\ \hline
IV. & (a) & spacelike line & & $(E_1, 0)$ \\ \cline{2-3} \cline{5-5}
& (b) & lightlike line & 0-degenerate & $(\frac12 (E_0+E_1), 0)$ \\
\cline{2-3} \cline{5-5}
& (c) & timelike line & & $(E_0, 0)$ \\ \cline{2-3} \cline{5-5}
& (d) & $(0,0)$ & & $(0,0)$ \\ \hline
\end{tabular} \end{center} \end{minipage}

\section{1-degenerate parallel surfaces in $\Rf$ }
\label{sec:1-degPar}

In the following only 1-degenerate surfaces in $\Rf$ will be
considered since only for this class the parallel surfaces are yet not
classified (cp. Sec.~\ref{sec:intro}).

\subsection{Second order frames}
\label{sec:soframes}

As we have seen before for a 1-degenerate surface (Typ II) in $\Rf$
there exists a frame $S\in F^1$ such that 
$h^3=E_2= \bigl( \begin{smallmatrix} 0&1\\ 1&0 \end{smallmatrix}
\bigr)$ and 
$h^4=\frac12 (E_0+E_1) = \bigl( \begin{smallmatrix} 1&0\\ 0&0
\end{smallmatrix} \bigr)$ resp. $\psi= \bigl( \begin{smallmatrix}
  \om^2 & \om^1\\ \om^1 &0 \end{smallmatrix} \bigr)$. We call such a
frame a {\em second order frame} and denote the set of all second
order frames on $U$ by $F^2\subset F^1$. We can determine the subgroup
$H^2\subset H^1$, which acts transitively and effectively on $F^2$, by
calculating which changes of frames leave the special form of $\psi$
invariant, and we obtain:

\begin{equation}
  \label{H2}
 H^2=\left\{ \left( \ba{c|c}
    \ba{c|c} \ba{cc} a & b\\ 0 & c\\ \ea   & 0\\ \hline 
      0&\ba{cc} ac&0 \\ 2ab & a^2\\ \ea  \\ \ea & 0 \\ \hline 
      0 & 1 \\ \ea \right) : ac \neq 0 \right\}.
\end{equation}

We have constructed a subbundle $\F^2$ of $\F^1$, $\F^2=\pi_U \colon
F^2 \to U$, which is a principal (right) $H^2$-bundle, the reduced
bundle obtained by reduction of the structure group $H^1$ of $\F^1$ to
$H^2$. We use the notation $S\in \F^2$ for a second order frame
field. The structure equations have the form
(cp.~\eqref{strucEqFone}): 

\begin{equation}
  \label{strucEqFtwo}
  d\spafi{v_1}{v_2}{\xi_1}{\xi_2}{x}=
\left( \ba{c|c}
% Achtung, ziemliche Bastelei!!
    \ba{c|c} \quad \iz\quad &
      \ba{cc} \om^2 & \om^1 \\ \om^1 & 0 \\ \ea \\ \hline 
        \quad \sigma\quad & \tau \\ \ea & 0 \\ \hline 
      \hspace{-1em}\om^1\ \; \om^2\ \,\;\quad 0\quad\,\; 0 & 0 \\ 
    \ea \right)
  \spafi{v_1}{v_2}{\xi_1}{\xi_2}{x}.
\end{equation}

\subsection{Parallel surfaces}
\label{sec:parallel}

An affine surface with transversal bundle $\sig$ is called {\em
parallel} if the second fundamental form $h= h^3 \xi_1 + h^4\xi_2$ is
parallel (cp. Sec.~\ref{sec:intro}), i.~e.
\begin{equation}
  \label{parallel}
  \nabla h=0 .
\end{equation}
By definition ($\nabla h = C^3 \xi_1 + C^4 \xi_2$, cp.~\cite{NV}) this
is equivalent to the vanishing of the cubic forms. In the following we
will use the abbreviation: $\iz_{ji}^k=\iz_j^k(X_i)$.

\begin{lem}\label{lem:parallel}
  If $(x,\sig)$ is a 1-degenerate parallel surface in $\Rf$, then we get for
  a second order frame field:
\begin{align}
\iz_{21}^1=0, \quad \iz_{22}^1=0, \\
\np_{X_1}\xi_1 &= (\iz_{11}^1 +\iz_{21}^2)\xi_1, \\ 
\np_{X_2}\xi_1 &= (\iz_{12}^1 +\iz_{22}^2)\xi_1, \\
\np_{X_1}\xi_2 &= 2\iz_{11}^2\xi_1 +2\iz_{11}^1\xi_2, \\
\np_{X_2}\xi_2 &= 2\iz_{12}^2\xi_1 +2\iz_{12}^1\xi_2.
\end{align}
\end{lem}
\begin{proof}
  This is a direct consequence of \eqref{parallel}, using \eqref{equ3}
  and $h_{22}=0$, $h_{12}=\xi_1$ and $h_{11}=\xi_2$ ($h_{ij}:=
  h(X_i,X_j)$).
\end{proof}

\begin{lem}
\label{ruled}
A 1-degenerate parallel surface in $\Rf$ is a ruled surface.
\end{lem}
\begin{proof} 
$D_{X_2} dx(X_2)= dx(\nabla_{X_2} X_2)+h(X_2,X_2)= \iz_{22}^2 X_2$.
\end{proof}

\begin{rem}
Every ruled surface $x(u,v)= \alpha(u)+ v \beta(u)$ in $\Rf$ is
$k$-degenerate ($k\in\{0,1\}$).
\end{rem}

\subsection{Further adaption of the frame and the parametrization}
\label{sec:adaption}

We know by now that a 1-degenerate parallel surface in $\Rf$ is a
ruled surface where $X_2$ gives the direction of the ruling. To
simplify the computations we would like to find a second order frame
field such that $\nabla_{X_2} X_2 =0$ and such that $\{X_1, X_2\}$ is
a Gauss-basis (i.~e. $0=[X_1, X_2]$).

\begin{lem}\label{lem:adapted}
  For a 1-degenerate parallel surface $(x,\sig)$ in $\Rf$ there exist
  a frame field $S=\{v_1, v_2, \xi_1,\xi_2, x\} \in \F^2$ and local
  coordinates $(u,v)$ such that $v_1=dx(\pu)$, $v_2=dx({\pv})$ and
  $D_\pv dx(\pv) =0$. Furthermore we can parametrize the surface by
  $x(u,v)= \alpha(u)+ v \beta(u)$.
\end{lem}
\begin{proof} 
If $S, \tS\in \F^2$ and $v_i=dx(X_i)$ resp. $\tv_i=dx(\tX_i)$,
then there exists $A\in H^2$ with $\tX_1= a X_1+ b X_2$, $\tX_2 = c
X_2$ and $ac\neq 0$ (cp.~\eqref{H2}). Thus
\begin{align*}
  0=\nabla_{\tX_2} \tX_2 &\iff X_2(\ln c)= -\ic222  \\
  0=[\tX_1,\tX_2] &\iff \begin{cases} X_2(\ln a)= -\ic121,\\
    X_1(\ln c)+\frac{b}{a}X_2(\ln c) -\frac1a X_2(b) = \ic122 -\ic212 .
  \end{cases}
\end{align*}
\end{proof}

We call a frame field $S=\{dx(\pu), dx({\pv}), \xi_1,\xi_2, x\}\in
\F^2$ with $D_\pv dx(\pv) =0$, where $(u,v)$ are local coordinates, an
{\em adapted frame field}. Let $x(u,v)= \al(u)+ v \be(u)$ be a
local parametrization of a 1-degenerate parallel surface in $\Rf$ and
$S=\{\al'+v\be', \be, \xi_1,\xi_2, x\}$ an adapted frame field. (For a
function $f(u)$ we write $f'=\pu f$.) By Lemma~\ref{lem:parallel}
and Lemma~\ref{lem:adapted} the structure equations have the following
form:
\begin{alignat}{2}
  \al''+v\be''=D_{\pu} dx(\pu) &= \ic111 (\al'+v\be') + \ic112 \be &+&
     \xi_2, \label{Guu}\\ 
  \be'=D_{\pv} dx(\pu) &= \ic212 \be &+&
     \xi_1, \label{Guv}\\ 
  D_{\pv} dx(\pv) &= 0, && \label{Gvv}\\
  D_{\pu} \xi_1 &= dx(-A_{\xi_1}(\pu)) &+& (\ic111 + \ic212)\xi_1,
     \label{Wu1}\\ 
  D_{\pv} \xi_1 &= dx(-A_{\xi_1}(\pv)), & & \label{Wv1}\\
  D_{\pu} \xi_2 &= dx(-A_{\xi_2}(\pu)) &+& 2\ic112 \xi_1 + 2\ic111 \xi_2,
     \label{Wu2}\\ 
  D_{\pv} \xi_2 &= dx(-A_{\xi_2}(\pv)) &+& 2\ic212
     \xi_1.\label{Wv2}
\end{alignat}
From \eqref{Guv} we get that $\xi_1= \be'- \ic212 \be$. Inserted in
\eqref{Wu1} this gives that $\be''-(\ic111+2\ic212)\be'$ must be
tangential. Since $dx(\pu)$, $dx(\pv)$ and $\xi_1$ are linear
independent, $\be'',\al', \be'$ and $\be$ must be linear dependent. We
get two cases:
\begin{align}
  \label{dalpha}
  {\rm I.}\quad& \al'(u)=k_1(u)\be(u)+ k_2(u)\be'(u)+ k_3(u)\be''(u),\\
  \label{ddbeta}
 {\rm II.}\quad& \be''\in \Span(\be,\be').
\end{align}
I. We assume that \eqref{dalpha} is true. We will investigate if it is
possible to reparametrize the surface such that $\al'=\be''$. If we
have coordinates $(\tu,\tv)\in \tU$ and a parametrization
$\tx(\tu,\tv)= \tal(\tu)+ \tv \tbe(\tu)$, then we can reparametrize
the surface by a local diffeomorphism $\phi\colon U \to \tU$,
$\phi(u,v)=(f(u), g(u)+ v h(u))$. We get $\tx \circ \phi(u,v)=: x(u,v)
=: \al(u) + v \be(u)$ with
\begin{align}
\al(u) &= \tal(f(u)) + g(u) \tbe(f(u)), \label{repal}\\
\be(u) &= h(u) \tbe(f(u)), \label{repbe} \\
h(u) f'(u) &\neq 0 \quad\forall u\in U. \label{diffeom}
\end{align}
Obviously the frame $S=(\al'+ v\be', \be,\xi_1, \xi_2)$ is an
adapted frame if\mbox{}f $\tS=(\tal'+ \tv \tbe', \tbe, \tilde{\xi}_1,
\tilde{\xi}_2)$ is an adapted frame. (If $J\phi$ is the
Jacobi matrix of $\phi$, then $S=B \tS$ with $B= \left( \ba{c|c}
    \ba{c|c} J\phi   & 0\\ \hline 0& \dots  \\ \ea & 0 \\ \hline 
      0 & 1 \\ \ea \right)\in H^2$.)
    
    Using \eqref{dalpha} resp. \eqref{repal}, \eqref{repbe} and
    \eqref{dalpha} for $\tal'$ ($:=\frac{\partial \tal}{\partial
      \tu}$), we obtain for the difference $\be''-\al'$ the following
    expression:
\begin{multline}
  \label{difference}
  \be'' -\al'=  -k_1\be -k_2 \be' +(1-k_3)\be''\\
  =(h''-g'-f' \tk_1)\tbe + (2h' f' +h f'' - gf'- f'
  \tk_2)\tbe' + (h (f')^2 -f' \tk_3)\tbe''.
\end{multline}
Thus $k_3\equiv 1$ iff $\tk_3 \circ f= h f'$. Hence we can find a
reparametrization such that $k_3\equiv 1$ (e.~g. $f=\text{id}$,
$h=\tk_3$) and it stays constant equal one if
we restrict to reparametrizations with $h=\frac{1}{f'}$, therefore 
\begin{equation}\label{restr1} h'=-\frac{f''}{(f')^2}. 
\end{equation} 
Now $k_2\equiv 0$ iff $\tk_2 \circ f= -\frac{f''}{(f')^2}-g$. Still we
  can find such a reparametrization (e.~g. $f=\text{id}$, $-\tk_2=g$)
  and $k_2$ stays constant equal zero if we restrict to
  reparametrizations with
\begin{equation}\label{restr2} g=-\frac{f''}{(f')^2} = h'.
\end{equation}
Finally $k_1\equiv 0$ iff $\tk_1 \circ f= \frac{1}{f'}(h''-g') \equiv
0$ (by \eqref{restr2}).

Since our investigations are of local nature, we have to consider two
subcases: either $\tk_1 \equiv 0$ or $\tk_1(\tu) \neq 0$ for all $\tu \in
\tU$. By \eqref{difference} this is equivalent to either $\be'' = \al'$
or $\be'' = \al' - (f')^2(\tk_1\circ f)\be$. In the second subcase we
still can choose $f$ such that $(f')^2(\tk_1\circ f)\equiv \pm 1$,
thus we have either
$$\al'= \be''\quad \text{or}\quad \al' = \be'' \pm \be.$$
Since $x(u,v)=\al(u)+ v \be(u)$, we obtain by integration (and if
necessary by an affine transformation applied to $x$) the following
two subcases:
\begin{align*}
{\rm I.1.}\quad &x(u,v)= \ga'(u)+ v \ga(u)\quad \text{or}\\
{\rm I.2.}\quad & x(u,v)= (\pm\ga(u) +\ga''(u))+ v \ga'(u).
\end{align*}

\noindent II. We assume that $\be''(u)\in \Span(be, be')(u) \quad\forall
u$. Therefore $\be$ is a plane curve, contained in the plane spanned
by $\be(u)$ and $\be'(u)$ for some $u$. We can reparametrize the
surface by a local diffeomorphism $\phi\colon U \to \tU$,
$\phi(u,v)=(f(u), v h(f(u)))$ (cp. the discussion in the first case)
such that $\be$ is part of an ellipse in $\Span(\be, \be')$ and, by
applying an affine transformation, such that $\be(u)=(\cos u, \sin u,
0,0)$, i.~e. $\be''=-\be$. We obtain
\begin{equation*}
{\rm II.}\quad x(u,v)=\al(u)+ v \be(u), \quad\be''=-\be.
\end{equation*}

To complete our investigations we need to compute for the three types
of ruled surfaces the corresponding transversal bundle $\sig$. This
can be done using the structure equations \eqref{Guu} - \eqref{Wv2}.
The computations are lengthy but straightforward. We give only a short
outline.

\noindent I. Let $\ga$ be a smooth function on an open subset of $\R$
such that $G:=\det(\ga, \ga',\ga'',\ga''')\neq 0$. We set
$$L:=\ln G \quad \text{and}\quad \ga^{(4)}=: L'\ga'''+ a\ga'' + b\ga'
+ c\ga.$$

\noindent 1. ($x(u,v)= \ga'(u)+ v \ga(u)$): If
we compute the Gauss equations \eqref{Guu} and \eqref{Guv} we obtain
\begin{align}
\xi_1&=\ga'-\ic212 \ga,\label{Fxi1}\\
\xi_2&=\ga'''+(v-\ic111)\ga'' -v\ic111\ga' -\ic112\ga. \label{Fxi2}
\end{align}
If we differentiate $\xi_2$ in direction of $u$ and evaluate the
Weingarten equation \eqref{Wu2}, we get:
\begin{equation*}
\ic111= \frac13(L'+v),\quad \ic112= \frac13(b-v a+v^2(L'+v)).
\end{equation*}
The Weingarten equation \eqref{Wu1} for $\pu\xi_1$ finally gives:
\begin{equation*}
\ic212=- \frac16(L'+4v).
\end{equation*}

\noindent 2. ($x(u,v)= (\eps\ga(u) +\ga''(u))+
v \ga'(u)$, $\eps =\pm 1$): We obtain by \eqref{Guu} and \eqref{Guv}
that 
\begin{align}
\xi_1&=\ga''-\ic212 \ga',\label{Sxi1}\\
\xi_2&=(L'+v-\ic111)\ga'''+(a+\eps-v\ic111)\ga'' +(b-\eps
\ic111-\ic112) \ga' +c \ga.\label{Sxi2}
\end{align}
If we differentiate $\xi_2$ in direction of $u$ and evaluate the
Weingarten equation \eqref{Wu2}, we get:
\begin{align*}
\ic111&= \frac13((\ln c +L)'+v),\\ 
\ic112&= \frac13\left\{b +a' -\eps L'-(a+\eps)(\ln c)' + v \left(-L''+
L'(\ln c)' -a-2\eps\right) + v^2 (\ln c +L)' + v^3\right\}. 
\end{align*}
The Weingarten equation \eqref{Wu1} for $\pu\xi_1$ finally gives:
\begin{equation*}
\ic212=- \frac16((\ln c +L)'+4v).
\end{equation*}

\noindent II. Let $\al, \be$ be smooth functions on an open subset of $\R$
such that $D:=\det(\al'',\al',\be',\be)\neq 0$. We set
$$L:=\ln D \quad \text{and}\quad \al'''=: L'\al''+ a\al' + b\be' +
c\be.$$
If we compute the Gauss equations \eqref{Guu} and \eqref{Guv}
we obtain
\begin{align}
\xi_1&=\be'-\ic212 \be,\label{Txi1}\\
\xi_2&=\al''-\ic111\al' -v\ic111\be' -(v+\ic112)\be. \label{Txi2}
\end{align}
If we differentiate $\xi_2$ in direction of $u$ and evaluate the
Weingarten equation \eqref{Wu2}, we get:
\begin{equation*}
\ic111= \frac13 L',\quad \ic112= \frac13(b-v a-v).
\end{equation*}
The Weingarten equation \eqref{Wu1} for $\pu\xi_1$ finally gives:
\begin{equation*}
\ic212=- \frac16 L'.
\end{equation*}

\begin{thm}
Every 1-degenerate parallel affine surface immersion $(x,\sig)$ in
$\Rf$ is a ruled surface and can be locally parametrized either by
\begin{itemize} 
\item[{\rm I.1.}] $x(u,v)= \ga'(u) + v \ga(u)$, and 
  \newline $\sig=\Span(\xi_1,\xi_2)$ is
  given by \eqref{Fxi1} and \eqref{Fxi2}, or
\item[{\rm I.2.}] $x(u,v)= (\eps \gamma(u) + \gamma''(u)) + v \gamma'(u)$,
  $\eps=\pm 1$, and 
  \newline $\sig=\Span(\xi_1,\xi_2)$ is given by \eqref{Sxi1} and
  \eqref{Sxi2},
\item[{\rm II.}] $x(u,v)=\al(u)+v \be(u)$, $\be''=-\be$, and
  \newline $\sig=\Span(\xi_1,\xi_2)$ is given by \eqref{Txi1} and
  \eqref{Txi2}. 
\end{itemize}
\end{thm}

%%%%%%%%%%%%%%  Bibliography %%%%%%%%%%%%%%%%%%

\end{document}